\documentclass[a4paper, 12pt]{article}
\usepackage{amsmath,amsthm, amsfonts}
\usepackage{amssymb}
\usepackage{rotating}
\newtheorem{theorem}{Theorem}[section]
\newtheorem{lemma}[theorem]{Lemma}
\newtheorem{proposition}[theorem]{Proposition}
\newtheorem{corollary}[theorem]{Corollary}

\newtheorem{problem}[theorem]{Problem}
\theoremstyle{definition}
\newtheorem{definition}[theorem]{Definition}
\newtheorem{remark}[theorem]{Remark}
\newtheorem{example}[theorem]{Example}

\newcommand{\Q}{\mathbb{Q}}
\newcommand{\C}{\mathbb{C}}
\newcommand{\X}{\mathcal{X}}
\newcommand{\N}{\mathbb{N}}
\newcommand{\Z}{\mathbb{Z}}
\newcommand{\F}{\mathbb{F}}

\newcommand{\NN}{\mathcal{N}}

\newcommand{\veps}{\varepsilon}
\newcommand{\GA}{\operatorname{GA}}
\newcommand{\Aff}{\operatorname{Aff}}
\newcommand{\TA}{\operatorname{TA}}
\newcommand{\TAM}{\operatorname{TA}}
\newcommand{\Bij}{\operatorname{Bij}}
\newcommand{\DA}{\operatorname{DA}}

\newcommand{\GL}{\operatorname{GL}}
\newcommand{\chr}{\operatorname{char}}
\newcommand{\FF}{\mathbb{F}}
\newcommand{\MA}{\operatorname{MA}}

\newcommand{\GTAM}{\operatorname{GTAM}}
\newcommand{\GLIN}{\operatorname{GLIN}}
\newcommand{\TLIN}{\operatorname{TLIN}}
\newcommand{\kar}{\operatorname{char}}
\newcommand{\XX}{\vec{X}}

\title{Polynomial automorphisms over finite fields: Mimicking non-tame and tame maps by the Derksen group.}
\author{\begin{tabular}{ll}
Stefan Maubach\footnote{Funded by Veni-grant of council for the
physical sciences, Netherlands Organisation for scientific
research (NWO).} & Roel Willems\footnote{Funded by Phd-grant of council for the
physical sciences, Netherlands Organisation for scientific
research (NWO).}\\
&\\
\small
Radboud University Nijmegen&\\
\small Postbus 9010, 6500 GL Nijmegen  &\\
\small The Netherlands&\\
\small s.maubach@math.ru.nl& \small r.willems@math.ru.nl\\
\end{tabular}}

\begin{document}

\maketitle

\begin{abstract}
If $F$ is a polynomial automorphism over a finite field $\F_q$ in dimension $n$, then it induces a bijection $\pi_{q^r}(F)$ of $(\F_{q^r})^n$ for every $r\in \N^*$. 
We say that $F$ can be `mimicked' by elements of a certain group of automorphisms $\mathcal{G}$ if there are $g_r\in \mathcal{G}$ such that $\pi_{q^r}(g_r)=\pi_{q^r}(F)$. 

In section \ref{fix} we show that the Nagata automorphism (and any other automorphism in three variables fixing one variable) can be mimicked by tame automorphisms. This on the one hand
removes the hope of showing that such an automorphism is non-tame by studying one bijection it induces, but on the other hand indicates that the bijections of tame automorphisms coincide with bijections of all automorphisms.

In section \ref{Derksen} we show that the whole tame group $\TA_n(\F_q)$ in turn can be mimicked by the Derksen subgroup $\DA_n(\F_q)$, which is the subgroup generated by affine maps and one particular element. It is known that if a field $k$ has characteristic zero, then $\TA_n(k)=\DA_n(k)$, which is not expected to be true for characteristic $p$.  

In section \ref{lin} we consider the subgroups $\GLIN_n(k)$ and $\GTAM_n(k)$ of the polynomial automorphism group (which are candidates to equal the entire automorphism group). An automorphism $F$ is called linearizable (in $\GA_n(k)$) if there is an automorphism 
$\varphi$ (in $\GA_n(k)$) such that $\varphi^{1}F\varphi$ is linear. What is generally 
denoted by ``the linearization conjecture'' is the question if $F\in \GA_n(k)$ such that 
$F^s=I$ implies that $F$ is linearizable. In this respect, it can be useful to have a method to distinguish some  non-linearizables from linearizables. In section \ref{lin} we explore such a method, and show that it can be of use only for automorphisms over $\F_2$: we show that $\GLIN_n(\F_q)=\GTAM_n(\F_q)$ if $q\not = 2$, and $\GLIN_n(\F_2)\subsetneq \GTAM_n(\F_2)$.
\end{abstract}

\section{Introduction}

Polynomial automorphisms are generally studied over $\C$, $\Q$, or any field of characteristic zero. Even if they are studied over commutative rings,
then it is often assumed that $\Q$ (or $\Z$) is a subset of the ring. The characteristic $p$ case is quite unexplored, though it is gaining some interest.

Denote by $\GA_n(k)$ the polynomial automorphism group in dimension $n$ over $k$, and $\TA_n(k)$ as the tame subgroup of $\GA_n(k)$ (a precise definition is given in the next section). Any element $F\in \GA_n(\F_q)$ (where $q=p^r$, $p$ a prime, and $\F_q$ denotes the finite field with $q$ elements)
induces a bijection of $\F_{q^m}^n$ for each $m\in \N^*$. This bijection we denote by $\pi_{q^m}(F)$.

One result from \cite{Maubach1} is that $\pi_q(\TA_n(F_q))$, $n\geq 2$, equals the set of all bijection of $(\F_q)^n$, except if 
$q=2^m, m\geq 2$, then any 
such bijection will be an even permutation of $\F_q^n$. 
This incited the search for automorphisms which were ``odd'', as such an example would 
immediately be non-tame, giving a very simple proof. Note that the proof of  $\TA_3(k)\not = \GA_3(k)$ by Umirbaev-Shestakov in \cite{UmSh1, UmSh2} 
is only in $\kar{k}=0$.
Unfortunately, all studied examples so far turned out to be even. 

The next idea was to study the group $\pi_{q^2}(\TA_n(\F_q))$, and find an automorphism $F\in \GA_n(F_q)$ such that
$\pi_{q^2}(F) \not \in \pi_{q^2}(\TA_n(\F_q))$. However, no such examples were found either.
In section \ref{fix} we give part of the reason that it is hard to find such an example: most possible counterexamples to being tame, fix one variable. 
It is of such maps, fixing one variable, that corollary \ref{MimickDim2} shows that they are always mimickable (as defined in the abstract).

This result may indicate that $\pi_{q^m}(\TA_n(\F_q))$ equals $\pi_{q^m}(\GA_n(\F_q))$. Such a result would not be without implication, as it would show that in any practical application over finite fields, it is enough to restrict to tame maps. In light of this question, it might even be useful to 
find an even smaller group that has this property. 
That is done in section \ref{Derksen}, where we introduce a (potentially) smaller group $\DA_n(\F_q)$ (the Derksen group) of which we prove that $\pi_{q^m}(\DA_n(\F_q))=
\pi_{q^m}(\TA_n(\F_q))$. It must be noted that if $\kar{k}=0$, then $\DA_n(k)=\TA_n(k)$, but for $\kar{k}=p$ such an equality is not expected to hold (see section \ref{TamGen}).

In the paper \cite{PolMau} the subgroups $\GLIN_n(k)$ and $\GTAM_n(k)$ of $\GA_n(k)$ are introduced. $\GLIN$ is the group generated by the set of automorphisms which are linear up to conjugation by an element of $\GA$, and $\GTAM$ is defined in a similar way.
In section \ref{lin} we show that these groups are equal, except if the field $k=\F_2$, in which case they differ. 

\section{Preliminaries and definitions}
Let $p$ be a prime, $q=p^r$, $\F _q$ the finite field with q elements. Let $n\geq 1$.
We are interested in the group $\GA_n(\F_q)$ of polynomial autmorphisms over $\F_q$.
Subgroups of $\GA_n(\F_q)$ are the group of linear automorphisms $\GL_n(\F_q)$ and the group of affine automorphisms $\Aff_n(\F_q)$.
Let us first fix a notation for some special elements of $\Aff_n(\F_q)$;
\begin{itemize}
 \item $T_{i,c} = (X_1,\ldots,X_{i-1},X_i+c,X_{i+1},\ldots ,X_n)$, 
 \item $S_{i,c} = (X_1,\ldots,X_{i-1},cX_i,X_{i+1},\ldots ,X_n)$ and 
 \item $R_{i,j} = (X_1,\ldots,X_{i-1},X_j,X_{i+1},\ldots ,X_{j-1},X_i,X_{j+1},X_n)$,
\end{itemize}
where $i,j\in \{1,\ldots ,n\}$ and $c \in \F_q^*$. Incidentally, the $T_{i,c},S_{i,c}, R_{i,j}$ generate $\Aff_n(\F_q)$.
Note that $R_{i,j}=R_{j,i}=R_{i,j}^{-1}$, $R_{i,j}T_{i,c}R_{i,j} = T_{j,c}$ and $R_{i,j}S_{i,c}R_{i,j} = S_{j,c}$.\\
Now if for $\alpha \in \N^{n-1}$ we write $X^{\alpha}=X_2^{\alpha_2}X_3^{\alpha_3}\cdots X_n^{\alpha_n}$ (note the omission of $X_1$),
then
we can define an elementary automorphism $E_{1,\alpha}=(X_1+X^{\alpha},X_2,\ldots,X_n)\in \GA_n(\F_q)$, and more general $E_{i,\alpha}=R_{1,i}E_{1,\alpha}R_{1,i}$.\\
Another subgroup of interest of $\GA_n(\F_q)$ is the group of tame automorphisms $\TA_n(\F_q)=\langle \Aff_n(\F_q), E_{1,\alpha};\alpha \in \N^{n-1}\rangle $ generated by the affine and elementary automorphisms.
An important result on this group is by Jung and van der Kulk \cite{Jung, Kulk}:
\begin{theorem}
For any field $k$, $\TA_2(k)=\GA_2(k)$ 
\end{theorem}
H. Derksen proved that 
\begin{theorem}\label{DT}
 Let $k$ be a field of characteristic zero and $n\geq 3$, then $\TA_n(k) = \langle \Aff_n(k), \veps\rangle$, where
$\veps = (X_1+X_2^2,X_2,\ldots,X_n)$.
\end{theorem}
For a proof, see \cite{Essenboek} page 95-96.
More recently Bodnarchuk in \cite{Bod} showed that $\TA_n(k)=\langle \Aff_n(k),F\rangle$, for $k$ a field of characteristic zero, where $F$ is any non-linear triangular automorphisms. As such, the choice of $\veps$ is quite arbitrary. However, over finite fields Bodnarchuk's result will for sure not hold, as shown in section \ref{TamGen}.
In Section \ref{Derksen} we will prove a similar (but weaker) result for $k$ a finite field.
Let us first define our version of the Derksen automorphism group, for an appropriate automorphism $\veps$:
\begin{definition}
\[ \DA_n(\F_q) = \langle \Aff_n(\F _q),\veps \rangle \subset \TA_n(\F_q),\]
where $\veps =E_{1,\alpha} = (X_1+X^{\alpha},X_2,X_3,\ldots ,X_n)$ where $\alpha:=(p-1,\ldots,p-1)$.
\end{definition}
We furthermore need the following map:
\begin{definition}
Let $\Bij(\F_q^n)$ be the group of bijections of $\F_q^n$, which is isomorphic to $Sym(q^{n})$. 
We can define
\[\pi_q:\GA_n(\F_q)\rightarrow \Bij(\F_q^n),\]
as the canonical map, associating to an automorphism its induced bijection on the space $\F_q^n$.
\end{definition}
Note that, since $\Bij(\F_q^n)\cong Sym(q^n)$, it follows that $\pi_q$ is a grouphomomorphism. In particular for $F,G\in \GA_n(\F_q)$ we have that $\pi_q(FG)=\pi_q(F)\pi_q(G)$. Furthermore since $\GA_n(\F_q)<\GA_n(\F_{q^m})$, we can talk about $\pi_{q^m}:\GA_n(\F_q)\rightarrow \Bij(\F_{q^m}^n)$.\\

\section{Automorphisms that fix a variable}
\label{fix}

{\bf Notations:} If $F\in \GA_n(k[Z])$ and $c\in k$, write $F_c\in \GA_n(k)$ for the restriction of $F$ to $Z=c$.
In case we already have a subscript $F=G_{\sigma}$, then $G_{\sigma,c}=F_c$.
If $F\in \GA_n(k)$, then by $(F,Z)\in \GA_{n+1}(k)$ (or any other appropriate variable in stead of $Z$) we denote the canonical map obtained by  
$F$ by adding one dimension. If $F\in \GA_n(k[Z])$ then we identify $F$ on $k[Z]$ with $(F,Z)$ on $k$ and denote both by $F$. (In fact, we think of 
$\GA_n(k[Z])$ as a subset of $\GA_{n+1}(k)$.)

\begin{definition}
Let $H$ be a subgroup of $\GA_n(\F_q)$.
If $F\in \GA_n(\F_q)$, then we say that $F$ is {\em $H$-mimickable}, or {\em can be mimicked by $H$-maps}, 
if for all $m\in \N$ there exists $G_m\in H$ such 
that $\pi _{q^m}(F)=\pi _{q^m}(G_m)$. 
\end{definition}

One could also state that the lattice $(G_i ~|~ i\in \N^*)$  is a ``$p$-adic approximation of $F$ by elements in $H$'' or something similar, but the name ``mimickable'' stuck more easily and we prefer it.

Before we state the actual theorem, let us give an example of an actual map that could be non-tame, but can be mimicked by tame ones.

\begin{example}
Let $f,g\in \F_q[Z]$. Define $\Delta:=gX+fY^2$, and 
\[ F:=(X-2fY\Delta-fg\Delta^2, Y+g\Delta) .\]
(The case that $g=Z,f=1$ is Nagata's automorphism.)
This map is in $\GA_2(\F_q[Z])$.
This map is tame over $\F_q(Z)$: \\
\[ F=(X-\frac{f}{g}Y^2,Y)(X,Y+g^2X)(X+\frac{f}{g}Y^2,Y). \]
In fact, it is tame over $\F_q[Z,g^{-1}]$. It is unlikely that $(F,Z)$ is tame over $\F_q$, however.
Fixing $m\in \N$, we will find a tame map $T\in \TA_2(\F_q[Z])$ for which $T_c=F_c$ for each $c\in \F_{q^m}$. 
First, define 
\[ G:=(X-fg^{q^m-2}Y^2,Y)(X,Y+g^2X)(X+fg^{q^m-2}Y^2,Y). \]
For any $c$ such that $g(c)\not =0$, $G_c=F_c$. 
If $g(c)=0$, then $G_c=(X,Y)$, but $F_c=(X-2f(c)Y^3,Y)$.
Now define $\rho:=(g^m-1)\in \F_q[Z]$. Note that $\rho(c)=0$ if $g(c)\not =0$, and $\rho(c)=1$ if $g(c)=0$.
Define
\[ H:= (X-2f(1-g^{q^m-1})Y^3,Y) .\]
Now if $q(c)\not = 0$, then $H_c=(X,Y)$. If $q(c)=0$, then $H_c=(X-2f(c)Y^3,Y)=F_c$.
Hence, we can take $T:=GH$ (or $HG$).
\end{example}

We will use the below remark several times in the proof of proposition \ref{Mimick}.\\

\begin{definition}
A map is called {\em strictly Jonqui\'ere} if it is Jonqui\'ere, sends 0 to 0, and has diagonal part identity. 
\end{definition}

\begin{remark}\label{remark}
 Let $F\in \TA_n(k(Z))$ be such that the affine part of $F$ is the identity. 
Then $F$ can be written as a product of strictly Jonqui\'ere maps and permutations.
\end{remark}

Proofs like the one below use arguments that can be called standard by those familiar with using the Jung-van der Kulk theorem.
Together with the fact that a precise proof is less insightful and involves a lot of bookkeeping, we decided to give a slightly rough sketch  proof:

\begin{proof} (rough sketch.)
The whole proof works since one can ``push'' elements which are both Jonqui\'ere and affine to one side, since if $E$ is Jonqui\'ere (or affine),
and $D$ is both, then there exists an $E'$ which is Jonqui\'ere (or affine), and $ED=DE'$. 
This argument is used to standardize many a decomposition. Here, we emphasize that the final decomposition is by no means of minimal length.\\ 
(1) First, using the Jung-van der Kulk theorem, we decompose $F=E_1A_1E_2A_2\cdots E_sA_s$ where each $E_i$ is Jonqui\'ere and 
each $A_i$ is affine. \\
(2) {\em We may assume that $E_i(0)=A_i(0)=0$ for all $1\leq i\leq s$.} For any pure translation part can be pushed to the left, and then we use the fact that $F(0)=0$. I.e. 
the $A_i$ are linear.\\
(3) {\em We may assume that $\det(A_i)=1$ and the $E_i$ are strictly Jonqui\'ere.} To realize this, one must notice that there exists a diagonal linear map $D_i$ satisfying $\det(D_i)=\det(A_i)$, and that we can do this by pushing diagonal linear maps to the left. The result follows 
since the determinant of the linear part of $F$ is 1, and hence the determinant of the Jacobian of $F$ is 1.\\
(4) {\em We may assume that each $A_i$ is either diagonal of determinant 1 or -1, or a permutation.}
Now we use Gaussian Elimination to write each $A_i=P_{i1}E_{i1}\cdots P_{it}E_{it}D_{i}$ as a composition of permutations $P_{ij}$,  strictly Jonqui\'ere (elementary linear) maps $E_{ij}$, and one  diagonal map $D_{i}$. We may assume that each $E_{ij}$ is in fact upper triangular, by  conjugating with a permutation. Note that the determinant of each $E_{ij}$ is 1, and of each $P_{ij}$ is 1 or -1, so the determinant of $D_i$ is 1 or -1.\\
(5) {\em We may assume that each $A_i$ is a permutation.} We have to replace the diagonal linear maps $D_i\in \GL_n(k(Z))$ which have determinant 1 or -1. First, write the diagonal linear map $D_i=D_{i1}\cdots D_{it}\tilde{D}_i$, where the $D_{ij}$ are diagonal linear of determinant 1  that have 1's on n-2 places, and at most two diagonal elements which are not 1, and $\tilde{D}_i$ has 1 on the diagonal except at one place, where it is 1 or -1. 
The following formula explains how to write a diagonal map $D_{it}$ as product of (linear) strictly Jonqui\'ere maps and permutations:
\[ 
\begin{array}{l}
\left( \begin{array}{cc}
f^{-1}& 0\\
0 & f\\
\end{array} \right)
=\\
\left( \begin{array}{cc}
1& f^{-1}\\
0 & 1\\
\end{array} \right)
\left( \begin{array}{cc}
0& 1\\
1 & 0\\
\end{array} \right)
\left( \begin{array}{cc}
1& 1-f\\
0 & 1\\
\end{array} \right)
\left( \begin{array}{cc}
0& 1\\
1 & 0\\
\end{array} \right)
\left( \begin{array}{cc}
1& -1\\
0 & 1\\
\end{array} \right)
\left( \begin{array}{cc}
0& 1\\
1 & 0\\
\end{array} \right)
\left( \begin{array}{cc}
1& 1-f^{-1}\\
 0& 1\\
\end{array} \right).
\end{array}\]
Finally, we give the reader the funny puzzle to  write $\tilde{D}$ as a product of linear strictly Jonqui\'ere maps and permutations, which indeed is always possible.
\end{proof}

The main result of this section is the following proposition:

\begin{proposition} \label{Mimick}
Assume that for some $m\in \N$ we have the following: 
If $F\in \TA_n(\F_q(Z))\cap \GA_n(\F_q[Z])$, and $F_c\in \TA_n(\F_q[c])$ for all $c\in{\F}_{q^m}$,
then   we can find $T\in \TA_n(\F_{q^m})$ such that $F_c=T_c$ for all $c\in \F_{q^m}$.
\end{proposition}

(Note that by theorem \ref{main} one can replace $\TA_n(\F_{q^m})$ by $\DA_n(\F_{q^m})$ even.)
This immediately yields the following important corollary:

\begin{corollary} \label{MimickDim2}
Let $F\in \GA_3(\F_q)$ fixing one variable. Then $F$ is tamely mimickable. 
\end{corollary}

In particular, it is impossible to show that the Nagata automorphism over a finite field $\F_p$  is non-tame by showing that
the bijection it defines over $\F_{p^m}$ is not an element of $\pi _{q^m}(\TA_3(\F_p))$.

\begin{proof} (of proposition \ref{Mimick})\\
{\bf (1)} We may assume that the affine part of $F$ is the identity, by, if necessary, composing with a suitable affine map.\\
{\bf (2)} Since $F\in \TA_n(\F_q(Z))$, we can use remark \ref{remark} and decompose $F$ into strictly Jonqui\'ere maps over $\F_q(Z)$ and permutations. 
Gathering all denominators which appear in this decomposition, we can assume that $F\in \TA_n(\F_q[Z,g(Z)^{-1}])$ for some
$g(Z)$. We may also assume that $g$ is radical, as $k[Z,g(Z)^{-1}]=k[Z,\textup{rad}(g(Z))^{-1}]$. \\
{\bf (3)} Using lemma \ref{OpenSetMimick} we find $G\in \TA_n(\F_{q}[Z])$ such that $F_{c}=G_c$ for all $c\in \F_{q^m}$ satisfying $g(c)\not =0$.
Replacing $F$ by $G^{-1}F$, we may:\\
{\bf (4)} Assume that $F_c=I_n$ for each $c$ such that $g(c) \not = 0$.\\
We decompose $g:=g_1g_2\cdots g_t$ into irreducible factors $g_i$.
For each $g_i$ we find a map $G_i$ as described in lemma \ref{ClosedSetMimick}.
Defining $T:=G_1G_2\cdots G_t$ we are done.
\end{proof}

\begin{lemma} \label{ClosedSetMimick}
Let $g(Z)\in \F_q[Z]$ be irreducible. 
Let $F\in \TA_n(\F_q[Z,g(Z)^{-1}])\cap \GA_n(\F_q[Z])$ be such that the affine part is the identity.
Assume that 
$F_c\in \TA_n(\F_q[c])$ for all $c\in{\F}_{q^m}$.
Then for any $m\in \N^*$ we find $G\in \TA_n(\F_q[Z])$ such that for $c\in \F_{q^m}$:\\
(1) $G_{c}=I_n$ if $g(c)\not =0$,\\
(2) $G_{c}=F_{c}$ if $g(c)=0$. 
\end{lemma}

\begin{proof}
We may assume that $m$ is such that $g$ factors completely into linear factors over $\F_{q^m}$ (for the result for divisors of $m$ is implied by the result for $m$).
Let $\alpha$ be a root of $g$. 
Consider $F_{\alpha}$, which by assumption and remark \ref{remark} 
can be written as a composition of strictly Jonqui\'ere maps $e_i$ and permutations $p_i$ : $F_{\alpha}=e_1p_1e_2\ldots e_sp_s$. 
Write the $e_i$ as $I_n+H_i$ where $H_i$ is strictly upper triangular.
We can even write $H_i=f_i(\alpha,\XX)$ where $f_i(Z,\XX) \in \F_q[Z,\XX]^n$ (and $\XX$ stands for $X_1,\ldots,X_n$).

Now we define $\rho:=1-g^{q^m-1}\in \F_q[Z]$ and $E_i:=I_n+ \rho f_i(Z, \XX)$. 
Note that all $E_i \in \TA_n(\F_q[Z])$.
We define $G:=E_1p_1E_2p_2\cdots E_sp_s$. Our claim is that this map acts as required.  Since for $c\in \F_{q^m}$ we have 
$\rho(c)=0$ if and only if $g(c)\not =0$, it follows that in that case
$G_c=I_n$. 
Since $E_{i,\alpha}=e_i$  by construction, we have $G_{\alpha}=F_{\alpha}$. 
Now let $\Phi$ be an element of the galois group $Gal(\F_{q^m}:\F_{q})$. The remaining question is if $G_{\Phi(\alpha)}=F_{\Phi(\alpha)}$. 
Note that if $P(\XX,Z)\in \F_q[\XX,Z]$ then $\Phi(P(\alpha,\XX))=P(\Phi(\alpha),\XX)$. This implies that if $F\in \GA_n(\F_q[Z])$, then 
$F_{\Phi(\alpha)}=\Phi(F_{\alpha})$. Thus $G_{\Phi(\alpha)}=\Phi(G_{\alpha})=\Phi(F_{\alpha})=F_{\Phi(\alpha)}$ and
we are done.
\end{proof}

\begin{lemma} \label{OpenSetMimick}
Let $g(Z)\in \F_q[Z]$ be irreducible. 
Let $F\in \TA_n(\F_q[Z,g(Z)^{-1}])\cap \GA_n(\F_q[Z])$ be such that the affine part is the identity.
Assume that 
$F_c\in \TA_n(\F_q[c])$ for all $c\in{\F}_{q^m}$.
Then for any $m\in \N^*$ we find $G\in \TA_n(\F_q[Z])$ such that for $c\in \F_{q^m}$:\\
$G_{c}=F_{c}$ if $g(c)\not =0$.\\
\end{lemma}

\begin{proof}
Using remark \ref{remark}, write $F=e_1a_1e_2a_2\cdots e_sa_s$ where the  $e_i$ are strictly Jonqui\'ere maps
and the $p_i$ are permutations. Write 
$e_i$ of the form $I_{n+1}+H_i$ where $H_i$ is strictly upper triangular.
Now we modify $H_i$ in the following way: replace each fraction $g^{-t}$ by $g^{t(q^m-2)}$, making new elements 
$\tilde{H}_i$ which are in $\MA_n(\F_q[Z])$. Write $E_i:=I_{n}+\tilde{H}_i$, and
define $G:=E_1p_1E_2p_2\cdots E_sp_s$.
Note that if $c\in \F_{q^m}$ and $g(c)\not =0$,
then $g^{-t}(c)=g^{t(q^m-2)}(c)$, thus also $E_{i,c}=e_{i,c}$.
In fact, the latter remark implies that $G_c=F_c$ for all $c\in \F_{q^m}$ such that $g(c)\not =0$. 
\end{proof}

\section{Generators of the Tame Automorphism Group}
\label{TamGen}

In characteristic zero Derksens theorem shows that one can generate the tame automorphism group by the affine maps and only one nonlinear map; Bodnarchuk's theorem states that for any nonlinear triangular map this is true. However:

\begin{remark} Bodnarchuk's theorem (and Derksen's theorem in its original form) is not true in characteristic $p$ in general: 
the simplest counterexample is $\pi_2(\Aff_3(\F_2), (X+Y^2,Y,Z))=
\pi_2(\Aff_3(\F_2))$ which consists of only even permutations, while $\pi_2(\TA_3(\F_2))$ consists of all bijections of $(\F_2)^3$.
\end{remark}

In fact, in characteristic $p$ it is not clear if there are finitely many automorphisms that one can add to the affine group to generate the whole tame automorphism group. We were able to find the following generating set $E$:

\begin{theorem}
$\TA_n(F_q)=\langle \Aff_n(F_q), E\rangle$,\\ 
where $E=\{(x_1+x_2^{k_2p-1}\cdots x_n^{k_np-1},x_2,\ldots , x_n)|1\leq k_2 \leq \ldots\leq k_n\}$.
\end{theorem}

The proof of this theorem is the topic of the current section.

\begin{proof}
 It suffices to show that $E_{1,v}\in \langle \Aff_n(F_q), E\rangle$ for all $v \in \N^{n-1}$.
We will proceed by induction to $v$, with respect to the standard lexicographic ordering on $\N^{n-1}$.
So fix  $v \in \N^{n-1}$ and let $k_2,\ldots ,k_n$ such that $(k_i-1)p\leq v_i\leq k_ip-1$.
By a conjugation with a suitable permutation we may assume that $v_2\leq v_3 \leq \ldots \leq v_n$ and $k_2\leq\ldots \leq k_n$.
Now from Lemma \ref{lem2a}, it follows that there exists a vector $\alpha = (\alpha_0,\ldots ,\alpha _{p-1})\in \F _p^p$ such that 
\[
 (\alpha_0,\ldots ,\alpha_{p-1})\left(\begin{array}{c}
(Y+0)^{k_np-1}\\
(Y+1)^{k_np-1}\\
\vdots\\
(Y+p-1)^{k_np-1}
\end{array}
\right)= Y^{((k_n-1)p+q)} + P(Y) 
\]
where $v_n=(k_n-1)p+q$ and $deg(P(Y))\leq (k_n-1)p-1$.\\
This means that if we let $F_i=S_{1,\alpha_i}T_{n,-i}E_{1,(k_2p-1,\ldots,k_np-1)}T_{n,i}S_{1,\alpha_i^{-1}}$, then
$F_0\circ\cdots \circ F_{p-1}=(X_1+X_2^{k_2p-1}\cdots X_{n-1}^{k_{n-1}p-1}\left(X_n^{v_n}+P(X_n)\right),X_2,\ldots , X_n)\in \langle \Aff_n(F_q), E\rangle$.\\
Now because $deg(P(y)\leq (k_{n}-1)p-1$  we have by induction that
\[(X_1+X_2^{k_2p-1}\cdots X_{n-1}^{k_{n-1}p-1}X_n^{v_n},X_2,\ldots , X_n)\in \langle \Aff_n(F_q), E\rangle\]
By repeating this procedure for $X_{n-1},\ldots ,X_2$, we get that
\[ E_{1,v}\in \langle \Aff_n(F_q), E\rangle\]
which proves our statement.
\end{proof}

\begin{lemma}\label{lem3}\label{lem1a}\label{lem2a}
Let $p$ be a prime, $q=p^r$ and $\F_q$ be the finite field with $q$ elements.
Let $\F_q[Y]$ be the ring in one variable $Y$, over $\F _q$. Let $k\in \N$, and finally, let $kp\leq l < kp+p$ where $l\in \N$.

Then there exists a vector $\alpha  =(\alpha_0,\ldots ,\alpha_{p-1})\in \F_p^p$ such that
\[
 (\alpha_0,\ldots ,\alpha_{p-1})\left(\begin{array}{c}
(Y+0)^{kp+p-1}\\
(Y+1)^{kp+p-1}\\
\vdots\\
(Y+p-1)^{kp+p-1}
\end{array}
\right)= Y^{l} + P(Y) 
\]
where  $deg(P(Y)) < kp$.
\end{lemma}
\begin{proof}
We will calculate modulo the $\F_q$-module $M$ of polynomials of degree $< kp$. 
First note that $(Y+i)^{kp+p-1} = \sum _{j=0}^{p-1} {kp+p-1 \choose j} i^jY^{kp+p-1-j} \mod{M}$. So\\
\[
\left( \begin{array}{c}
(Y+0)^{kp+p-1}\\
\vdots\\
(Y+p-1)^{kp+p-1}
\end{array}\right)  = \left({kp+p-1  \choose j}i^j\right)_{0\leq i,j\leq p-1} 
\left( \begin{array}{c}
Y^{kp+p-1}\\
\vdots\\
Y^{kp}
\end{array}\right) \mod{M}.
\]
Now
\[
 \left({kp+p-1 \choose j}i^j\right)_{0\leq i,j\leq p-1} = \left( \begin{array}{ccc}
                                                               {kp+p-1 \choose 0} & & \emptyset\\
                                                                & \ddots & \\
                                                               \emptyset && {kp+p-1 \choose p-1}
                                                              \end{array}\right)
\left(i^j\right)_{0\leq i,j\leq p-1}.
\]
Because ${kp+p-1 \choose j}\not= 0 \mod{p}$ for $0\leq j\leq p-1$, and 
$\left( i^j\right)_{0\leq i,j\leq p-1}$ is a vandermonde-matrix, and  invertible,
 it follows that this
is an invertible matrix. We can take $\alpha$ to be the $l$-th column of the inverse of this matrix.
\end{proof}

\section{Derksen Automorphisms as bijections}\label{Derksen}

In this section we will prove the following weaker version of Derksen's Theorem \ref{DT}:

\begin{theorem}\label{main}
Let $q=p^r$ and $\F _q$ the finite field with $q$ elements, and $n\geq 3$. Then
  \[\pi_{q^m}(\TA_n(\F _q))=\pi_{q^m}(\DA_n(\F_q))\]
 \end{theorem}

In other words, the Derksen group mimicks the tame automorphism group.
It is  quite unplausible to expect that $\DA_n(\F_q)=\TA_n(\F_q)$ - though unproven. 
The focus of this section is to prove the above theorem. Unfortunately, certain parts of the proof are unavoidably quite technical.

Before we can prove this we need some intermediate results, which we will prove first. 
The next proposition is our main tool in the proof of Theorem \ref{main}.

\begin{proposition}\label{prop0}
 Let $\X = \F _{q^m}^n$ and $\tilde{\X}=\X\setminus \{u\in \F_{q^m}^n|u_n=0\}$, then define the bijection $\psi:\tilde{\X}\rightarrow \tilde{\X}$ by $\psi (u_1,\ldots ,u_n) = (u_1u_n^{-1},u_2u_n,u_3,\ldots, u_n)$. Then there exist a tame automorphism $T_m\in \DA_n(\F _q)$ such that $\pi_{q^m}(T_m)|_{\tilde{\X}} = \psi$. 
\end{proposition}
\begin{proof} 
First we need some elements from $\Aff_n(\F _q)$, define
$G=E_{1,(1,0,\ldots,0)}$ and $H=R_{1,2}$.
There are two cases
\begin{itemize}
 \item $\chr (\F_q)=2$, or
 \item $\chr (\F _q)\not= 2$.
\end{itemize}
%If $\chr (\F _q) = 2$ then define $A=E_{1,(1,0,\ldots,0,1)}\in D_n(\F_q)$ (this follows from lemma 3).
%Now let $T_m=H(G((AH)^{2^{2m-1}})HGH)^2H = HG(AHAH)^{2m-1}HGHG(AHAH)^{2m-1}HG=
%(x_1x_n^{4^m-2}+x_2(1+\sum_{i=0}^{2m-1}x_2^{4^m-2^i-1}+x_n^{4m-1}),x_1(1+\sum_{i=0}^{2m-1}x_2^{4^m-2^i-1}+x_n^{4^m-1}) + %x_2(\sum_{i=1}^{2m} x_n^{4^m-2^i} + x_n^{4^m}),x_3,\ldots,x_n)$, proof by induction:\\
%First of all note that $B=AHAH = (x_1+x_2x_n+x_1x_n^2,x_2+x_1x_n,x_3,\ldots,x_n)$ and $F=HG=(x_2,x_1+x_2,x_3,\ldots ,x_n)$.
%$T_1 = H(GAHAHHGH)^2H=HGAHAHHGHGAHAHHG=FBFFBF=(FBF)^2$, and more general $T_m=(FB^{2m-1}F)^2$. Zit beetje vast hier..........\\
%\ \\
%\ \\
%then $\pi_{q^m}(T_m)_{|\tilde{X}}=\psi$.\\
%\ \\
%\ \\
\ \\
First the easier case where $\chr (\F _q)\not= 2$. It follows that for every $m\geq 1$ we can define $A_m=E_{1,(1,0,\ldots,0,q^m-2)}=(X_1+X_2X_n^{q^m-2},X_2,\ldots ,X_n)$, $B=E_{1,(1,0,\ldots,0)}\circ S_{1,-1}\circ E_{1,(1,0\ldots,0,1)}\circ S_{1,-1}=(X_1+X_2-X_2X_n,X_2,\ldots X_n)$, and $C_m=E_{1,(1,0,\ldots,0)}\circ S_{1,-1}\circ E_{1,(1,0,\ldots,0,q^m-2)}\circ S_{1,-1}=(X_1+X_2-X_2X_n^{q^m-2},X_2,\ldots,X_n)$.
Note that we need $-1(\not=0,1)\in \F_q$, which is the reason why this will not work for $\chr (\F_{q})=2$. 
From Corollary \ref{cor1} it follows that $A_m, B, C_m\in \DA _n(\F_q)$.
So it follows that $T_m=A_mHBHG^{-1}HC_mH \in \DA _n(\F_q)$.
One verifies that $T_m=(2X_1X_n^{q^m-2}-X_1X_n^{2q^m-3}+X_2X_n^{q^m-1}-X_2,X_1-X_1X_n^{q^m-1}+X_2X_n,X_3,\ldots,X_n)$.
Now for $t \in \F_{q^m}^*$ Fermat's little theorem states that $t^{q^m-1}=1$, so it follows that for $u\in \tilde{X}$ we have that $T_m(u) = (u_1u_n^{-1},u_2u_n,u_3,\ldots,u_n)$, which shows that $\pi_{q^m}(T_m)_{|\tilde{\X}}=\psi$.

If $\chr (\F_q)=2$, then let $A=E_{1,(1,0,\ldots,0,1)}=(X_1+X_2X_n,X_2,\ldots,X_n)\in \DA _n(\F_q)$ (this follows from Corollary\ref{cor1}).
Define $F = HG$ and $B_m=(AH)^{k_m}$, where $k_m = 2^{2m-1}$. We claim that $T_m=(FB_mF)^2$ satisfies all our requirements.\\
Write $h_m = \sum_{j=1}^{2m-1}X_n^{k_m-2^j}$.
From Lemma \ref{proofBm} we have that
 \[
\begin{array}{rl}
B_m  = & (X_1(X_n^{k_m}+h_m) + X_2X_n^{k_m-1},X_1X_n^{k_m-1} + X_2(h_m),X_3,\ldots , X_n)
\end{array}
\]\\
Now write $\tilde{F}=(X_2,X_1+X_2)$, so $F=(\tilde{F},X_3,\ldots ,X_n)$ and $FB_mF=(\tilde{F}\tilde{B}_m\tilde{F},X_3,\ldots ,X_n)$.\\
\[
 \begin{array}{rl}
  \tilde{F}\tilde{B}_m\tilde{F} = & \left( \begin{array}{cc}
0 &1\\
1&1
\end{array}\right)
\left( \begin{array}{cc}
X_n^{k_m}+h_m & X_n^{k_m-1}\\
X_n^{k_m-1} & h_m
\end{array}\right)
\left( \begin{array}{cc}
0&1\\
1&1
\end{array}\right)
\left( \begin{array}{c}
X_1\\X_2
\end{array}\right)\\
%=&
%\left( \begin{array}{cc}
%0 &1\\
%1&1
%\end{array}\right)
%\left( \begin{array}{cc}
%x_n^{k_m-1} & x_n^{k_m}+x_n^{k_m-1}+h_m\\
%h_m & x_n^{k_m-1} + h_m
%\end{array}\right)
% \left( \begin{array}{c}
%x_1\\x_2
%\end{array}\right)\\
=&
\left( \begin{array}{cc}
h_m & X_n^{k_m-1}+h_m\\
X_n^{k_m-1}+h_m & X_n^{k_m}
\end{array}\right)
 \left( \begin{array}{c}
X_1\\X_2
\end{array}\right)
\end{array}
\]
hence $T_m=(FB_mF)^2 = ((\tilde{F}\tilde{B}_m\tilde{F})^2,X_3,\ldots ,X_n)$. So
{\small{
\[
 \begin{array}{rl}
&(\tilde{F}\tilde{B}_m\tilde{F})^2  \\
=&\left( \begin{array}{cc}
h_m & X_n^{k_m-1}+h_m\\
X_n^{k_m-1}+h_m & X_n^{k_m}
 \end{array}\right)^2
\left( \begin{array}{c}
 X_1\\ X_2
 \end{array}\right)\\
= & 
\left( \begin{array}{cc}
 X_n^{2(k_m-1)}& X_n^{2k_m-1}+(X_n^{k_m} + X_n^{k_m-1}+h_m)h_m\\
X_n^{2k_m-1}+(X_n^{k_m} + X_n^{k_m-1}+h_m)h_m & X_n^{2k_m} + X_n^{2(k_m-1)} +h_m^2
 \end{array}\right)
\left( \begin{array}{c}
 X_1\\ X_2
 \end{array}\right)\\
=:&M\left(\begin{array}{c}
          X_1\\ X_2
          \end{array}\right)
\end{array}
\]}}
All we still have to show now is that $\pi_{2^m}(T_m)|_{\tilde{\X}} = \psi$.
For this we have to show that for $u\in \F_{2^m}^*$, we have that 
\begin{itemize}
 \item[1)] $u^{2(k_m-1)}=u^{-1}$.
\item[2)] $u^{2k_m-1}+(u^{k_m}+u^{k_m-1}+h_m(u))h_m(u)=0$.
\item[3)] $u^{2k_m}+ u^{2(k_m-1)}+h_m^2(u)=u$.
\end{itemize}
1) follows from Lemma \ref{lem123} i).
From Lemma \ref{lem123} ii) it follows that $h_m(u)=u^{2^{2m-1}-1}$, so 2) becomes
\[
 \begin{array}{rl}
 & u^{2k_m-1}+(u^{k_m}+u^{k_m-1}+h_m(u))h_m(u)\\
 =&u^{2^{2m}-1}+(u^{2^{2m-1}}+u^{2^{2m-1}-1}+u^{2^{2m-1}-1})u^{2^{2m-1}-1}\\
=&1+u^{2^{2m-1}}u^{2^{2m-1}-1}\\
=&1+u^{2^{2m}-1}\\
=&1+1=0
 \end{array}
\]
Finally 3)
\[
 \begin{array}{rl}
&  u^{2k_m}+ u^{2(k_m-1)}+h_m^2(u)\\
=& u^{2^{2m}} + u^{2^{2m}-2}+u^{2(2^{2m-1}-1)}\\
=& u +u^{-1}+u^{-1}\\
=& u
 \end{array}
\]
This proves that for the above matrix $M$, we have 
\[
 \begin{array}{rl}
  \pi_{2^m}(M)|_{\tilde{X}} =& \left( \begin{array}{cc}
X_n^{-1}&0\\
0&X_n
\end{array}\right)|_{\tilde{X}}
 \end{array}
\]
hence $\pi_{q^m}(T_m)|_{\tilde{\X}}=\psi$.
%\[
 %\begin{array}{rl}
%&  \pi_{2^m}(T_m)_{|\tilde{X}}\\
%=&( u_1u_n^{2(k_m-1)}+u_2(u_n^{2k_m-1}+(u_n^{k_m} + u_n^{k_m-1}+h_m(u_n))h_m(u_n)),u_1(u_n^{2k_m-1}+(u_n^{k_m} + %u_n^{k_m-1}+h_m(u_n))h_m(u_n))+u_2(u_n^{2k_m}+ u_n^{2(k_m-1)}+h_m^2(u_n)),u_3,\ldots ,u_n)\\
%=&(u_1u_n^{-1}+u_2\cdot0,u_1\cdot 0+ u_2u_n,u_3,\ldots ,u_n)\\
%=&(u_1u_n^{-1},u_2u_n,u_3,\ldots ,u_n)\\
%=&\psi
%\end{array}
%\]
%If $\chr (\F _q)\not= 2$, then for every $m\geq 1$ we can define $A_m=E_{1,(1,0,\ldots,0,q^m-2)}=(x_1+x_2x_n^{q^m-2},x_2,\ldots %,x_n)$, $B=E_{1,(1,0,\ldots,0)}\circ S_{1,-1}\circ E_{1,(1,0\ldots,0,1)}\circ S_{1,-1}=(x_1+x_2-x_2x_n,x_2,\ldots x_n)$, and %$C_m=E_{1,(1,0,\ldots,0)}\circ S_{1,-1}\circ E_{1,(1,0,\ldots,0,q^m-2)}\circ S_{1,-1}=(x_1+x_2-x_2x_n^{q^m-2},x_2,\ldots,x_n)$.
%Note that we need $-1\in \F_q$, which is the reason why this will not work for $\chr (\F_{q})=2$. 
%From Corollary \ref{cor1} it follows that $A_m, B, C_m\in \DA _n(\F_q)$.
%So it follows that $T_m=A_mHBHG^{-1}HC_mH \in \DA _n(\F_q)$.
%Using for instance a computer one easily verifies that %$T_m=(2x_1x_n^{q^m-2}-x_1x_n^{2q^m-3}+x_2x_n^{q^m-1}-x_2,x_1-x_1x_n^{q^m-1}+x_2x_n,x_3,\ldots,x_n)$.\\
%Now for $t \in \F_{q^m}^*$ Fermat's little theorem states that $t^{q^m-1}=1$, so it follows that for $u\in \tilde{X}$ we have that %$T_m(u) = (u_1u_n^{-1},u_2u_n,u_3,\ldots,u_n)$. Which shows that $\pi_{q^m}(T_m)_{|\tilde{X}}=\psi$.\end{proof}
\end{proof}

\begin{lemma}\label{proofBm}
If $\chr(\F_q)=2$. Let $A,F,G,H,B_m,k_m,h_m$ as in Proposition \ref{prop0}.
We will prove that 
\[
\begin{array}{rl}
B_m  = & (X_1(X_n^{k_m}+h_m) + X_2X_n^{k_m-1},X_1X_n^{k_m-1} + X_2(h_m),X_3,\ldots , X_n)
\end{array}
\]
\end{lemma}
\begin{proof} 
One can verify by an elementary computation that $B_1=(AH)^4 = (X_1+X_1X_n^2+X_2X_n,X_1X_n+X_2,X_3,\ldots,X_n)$.
Now notice that $k_{m+1}=2^{2(m+1)-1} =2^{2m-1+2}=2^2*2^{2m-1}=4k_m$, so $B_{m+1}=B_m^4$.
Write $B_m = (\tilde{B}_m,X_3,\ldots , X_n)$, where $\tilde{B}_m=(X_1(X_1^{k_m}+h_m)+X_2X_n^{k_m-1},X_1X_n^{k_m-1}+X_2h_m)$ and notice that $B_{m+1} = (\tilde{B}_{m+1},X_3,\ldots , X_n) = (\tilde{B}_{m}^4,X_3,\ldots , X_n) = B_m^4$. 
 so we have to verify that $\tilde{B}_{m}^4$ equals $\tilde{B}_{m+1}$
\[
\begin{array}{rl}
\tilde{B_m}^4 = & \left(
\begin{array}{cc}
X_n^{k_m}+h_m &  X_n^{k_m-1}\\
X_n^{k_m-1} & h_m
\end{array}\right)^4 
\left( \begin{array}{c} X_1 \\ X_2\end{array}\right) \\
% = & \left(
%\begin{array}{cc} 
%x_n^{4k_m}+h_m^4+ x_n^{2k_m}x_n^{2(k_m-1)}+x_n^{4(k_m-1)}&  x_n^{3k_m}x_n^{k_m-1}\\
%x_n^{3k_m}x_n^{k_m-1} & h_m^4+ x_n^{2k_m}x_n^{2(k_m-1)}+x_n^{4(k_m-1)}
%\end{array}\right)
% \left( \begin{array}{c} x_1 \\ x_2\end{array}\right) \\
% = & \left(
%\begin{array}{cc} 
%x_n^{4k_m}+h_m^4+ x_n^{4k_m-2}+x_n^{4k_m-4}&  x_n^{4k_m-1}\\
%x_n^{4k_m-1} & h_m^4+ x_n^{4k_m-2}+x_n^{4k_m-4}
%\end{array}\right) 
%\left( \begin{array}{c} x_1 \\ x_2\end{array}\right)\\
 = & \left(
\begin{array}{cc} 
X_n^{k_{m+1}}+h_m^4+ X_n^{k_{m+1}-2}+X_n^{k_{m+1}-4}&  X_n^{k_{m+1}-1}\\
X_n^{k_{m+1}-1} & h_m^4+ X_n^{k_{m+1}-2}+X_n^{k_{m+1}-4}
\end{array}\right) 
\left( \begin{array}{c} X_1 \\ X_2\end{array}\right)
\end{array}
\]
%recall that $\chr (\F_q) = 2$ and that $4k_m=k_{m+1}$.\\ 
We leave it to the reader to verify that $h_m^4+X_n^{k_{m+1}-2}+X_n^{k_{m+1}-4} = h_{m+1}$.
%\[
% \begin{array}{rl}
%  h_m^4+x_n^{k_{m+1}-2}+x_n^{k_{m+1}-4} =& (\sum_{j=1}^{2m-1}x_n^{k_m-2^j})^4 + x_n^{k_{m+1}-2} + x_n^{k_{m+1}-4}\\
% =& \sum_{j=1}^{2m-1}x_n^{4(k_m-2^j)}+x_n^{k_{m+1}-2} + x_n^{k_{m+1}-4}\\
%=& \sum_{j=1}^{2m-1}x_n^{4k_m-4*2^j} +x_n^{k_{m+1}-2} + x_n^{k_{m+1}-4}\\
%=& \sum_{j=1}^{2m-1}x_n^{k_{m+1}-2^{j+2}} +x_n^{k_{m+1}-2} + x_n^{k_{m+1}-4}\\
%=& \sum_{j=3}^{2m+1}x_n^{k_{m+1}-2^j} +x_n^{k_{m+1}-2^1} + x_n^{k_{m+1}-2^2}\\
%=& \sum_{j=1}^{2(m+1)-1}x_n^{k_{m+1}-2^j}\\
%=&h_{m+1}
% \end{array}
%\]
\end{proof}

\begin{lemma}\label{lem123}
 Let $q=2^m$, for some $m\geq 1$, let $u \in \F_q^*$ and let $h_m(X)=\sum_{j=1}^{2m-1}X^{2^{2m-1}-2^j}$, then the following statements are true;
\begin{itemize}
 \item[i)] $u ^{2^{2m}}=u$.
 \item[ii)] $h_m(u)=u^{2^{2m-1}-1}$.
\end{itemize}
\end{lemma}
\begin{proof}
Define $\varphi :\F _q^*\rightarrow \F_q^*$ as $\varphi (u)=u^2$, the Frobenius-automorphism.
First of all notice that $\varphi ^m(u)=u^{2^m}=u$. Now i) readily follows: $u ^{2^{2m}}=\varphi^{2m}(u)=\varphi^m(\varphi^m(u))=\varphi^m(u)=u$. 
For ii), define $v = u^{-1}$ and note that $u^{2^{2m}-1}=1$, then we have that
\[
 \begin{array}{rl}
h_m(u)=&\sum_{j=1}^{2m-1}u^{2^{2m-1}-2^j}\\
=&\sum_{j=1}^{2m-1}u^{2^{2m-1}}u^{-2^j}\\
=&u^{2^{2m-1}}\sum_{j=1}^{2m-1}v^{2^j}\\
=&u^{2^{2m-1}}\sum_{j=1}^{2m-1}v^{2^j} + 2u^{2^{2m-1}}v^{2^{2m}}\\
=&u^{2^{2m-1}}\sum_{j=1}^{2m}v^{2^j} + u^{2^{2m-1}}v^{2^{2m}}\\
=&u^{2^{2m-1}}\sum_{j=1}^{2m}\varphi^j(v) + u^{2^{2m-1}}v\\
=&u^{2^{2m-1}}(\sum_{j=1}^{m}\varphi^j(v) +\sum_{i=1}^{m}\varphi^{m+i}(v))+ u^{2^{2m-1}}u^{-1}\\
=&u^{2^{2m-1}}(\sum_{j=1}^{m}\varphi^j(v) +\sum_{i=1}^{m}\varphi^{i}(\varphi^m(v)))+ u^{2^{2m-1}-1}\\
=&u^{2^{2m-1}}(\sum_{j=1}^{m}\varphi^j(v) +\sum_{i=1}^{m}\varphi^{i}(v))+ u^{2^{2m-1}-1}\\
=&u^{2^{2m-1}}(2\sum_{j=1}^{m}\varphi^j(v))+ u^{2^{2m-1}-1}\\
=&u^{2^{2m-1}-1}
 \end{array}
\]
\end{proof}

\begin{lemma}\label{lem2}
 Let $\alpha=(1,\alpha_3,\ldots,\alpha_n)\in \N^{n-1}$ and $\beta=(0,\beta_3,\ldots,\beta_n)\in \N^{n-1}$, then 
\[
 [E_{1,\alpha},E_{2,\beta}]=E_{1,(0,\alpha_3+\beta_3,\ldots,\alpha_n+\beta_n)}
\]
\end{lemma}
\begin{proof}
Trivial.\end{proof}

\begin{lemma}\label{lem1}
 Let $\alpha=(\alpha_2,\ldots ,\alpha_n)\in \{0,\ldots ,p-1\}^{n-1}$, then $E_{1,\alpha}\in \DA_n(\F_q)$.
\end{lemma}
\begin{proof}
It follows from lemma \ref{lem3}, that for $l \in \{0,\ldots ,p-1\}$, there exist a vector $\beta _l=(\beta _l^0,\beta _l^1,\ldots ,\beta _l^{p-1})\in \F _p^p$, such that \[
 (\beta_l^0,\ldots ,\beta_l^{p-1})\left(\begin{array}{c}
(Y+0)^{p-1}\\
(Y+1)^{p-1}\\
\vdots\\
(Y+p-1)^{p-1}
\end{array}
\right)= Y^l
\]
In particular there exist such a vector for $l=\alpha_2$.
For $0\leq t\leq p-1$ define $c_t=\beta_{k_2}^t$ and let $F_{2,t}=S_{1,c_t}\circ T_{2,-t}\circ \veps \circ T_{2,t}\circ S_{1,c_t^{-1}} = (X_1+c_t(X_2+t)^{p-1}X_3^{p-1}\cdots X_n^{p-1},X_2,\ldots ,X_n)$.\\
So 
\[ \begin{array}{rl}
G_{\alpha_2}=&F_{2,0}\circ F_{2,1}\circ \cdots \circ F_{2,p-1}\\
 = &(X_1+(\sum_{i=0}^{p-1}c_t^i(X_2+i)^{p-1})X_3^{p-1}\cdots X_n^{p-1},X_2,\ldots ,X_n) \\
=& (X_1+X_2^{\alpha_2}X_3^{p-1}\cdots X_n^{p-1},X_2,\ldots ,X_n). 
\end{array} \]
Now repeating this for $l=\alpha_3$ with $T_{3,t}$ and $G_{\alpha_2}$ instead of $\veps$ and so on, gives us the required result.\end{proof}

\begin{proposition}\label{prop2}
 Let $\alpha=(0,\alpha_3,\ldots,\alpha_n)\in \N^{n-1}$, then $E_{1,\alpha}\in \DA_n(\F _q)$.
\end{proposition}
\begin{proof}
If $\alpha_i\leq p-1$ for all $3\leq i\leq n$, then the result follows from Lemma \ref{lem1}.
In particular it follows that $E_{1,(1,p-1,\ldots ,p-1)}\in \DA_n(\F _q)$.\\
Now we can use Lemma \ref{lem2} to construct $E_{1,\alpha}$ by induction.\end{proof}

\begin{corollary}\label{cor1}
 Let $\alpha \in \N^{n-1}$ with $\alpha _j = 0$ for some $j \in \{2,\ldots , n\}$, then $E_{i,\alpha}\in \DA_n(\F _q)$ if $i\not =j$.
\end{corollary}
\begin{proof}
Follows from Proposition \ref{prop2} and the fact that $R_{2,i}\in \DA_n(\F_q)$.\end{proof}

\begin{proposition}\label{prop1}
 Let $q,\ n,\ m, \X,\tilde{\X}, \psi$ and $T_m$ as in Proposition \ref{prop0}.
Let $\alpha = (\alpha_2,\ldots ,\alpha _n)\in \N^{n-1}$.
Now $\pi_{q^m}(T_m^{-1}\circ E_{1,\alpha}\circ T_m) = 
\pi_{q^m}(E_{1,\beta})$, where $\beta = (\alpha_2,\ldots ,\alpha_{n-1},\alpha_n+\alpha_2+1)$.
\end{proposition}
\begin{proof}
For $u\in \X\backslash \tilde{\X}$, $\pi_{q^m}(E_{1,\alpha})=\pi_{q^m}((x_1,\ldots ,x_n))=\pi _{q^m}(E_{1,\beta})$, so the statement is clearly true in $\X\backslash \tilde{\X}$. Now if we restrict ourselves to $\tilde{\X}$ we have that 
\[
 \begin{array}{rl}
   & \pi_{q^m}(T_m^{-1}\circ E_{1,\alpha}\circ T_m)\\
 = & \pi_{q^m}(T_m^{-1})\pi_{q^m}(E_{1,\alpha})\pi_{q^m}(T)\\
 = & \psi^{-1} \pi_{q^m}(E_{1,\alpha})\psi\\
 = & \pi_{q^m}((X_1X_n,X_2X_n^{-1},X_3,\ldots ,X_n))\pi_{q^m}(E_{1,\alpha})\pi_{q^m}((X_1X_n^{-1},X_2X_n,X_3,\ldots ,X_n))\\
 = &  \pi_{q^m}((X_1X_n,X_2X_n^{-1},X_3,\ldots ,X_n))\cdot\\
& ~~~~~\pi_{q^m}((X_1X_n^{-1}+X_2^{\alpha_2}\cdots X_{n-1}^{\alpha_{n-1}}X_n^{\alpha_2+\alpha_n},X_2X_n,X_3,\ldots,X_n))\\
 = & \pi_{q^m}(X_1+X_2^{\alpha_2}\cdots X_{n-1}^{\alpha_{n-1}}X_n^{\alpha_2+\alpha_n+1},X_2,X_3,\ldots,X_n)\\
 = & \pi_{q^m}(E_{1,\beta}) 
 \end{array}
\]
\end{proof}
\begin{lemma}\label{lemd}
 Let $a,b,m \in \Z$. If $Gcd(a,b,m)=d$, then there exists a $t\in \Z$ such that $Gcd(a+tb,m)=d$.
\end{lemma}
\begin{proof}
First assume $d=1$.
Define 
\[
 t=\prod_{\begin{tabular}{c}$p\ prime$\\ $p|m$\\ $p\not{|}a$\end{tabular}} p,
\]
then we show that $Gcd(a+tb,m)=1$: Let $p$ be a prime such that $p|Gcd(a+tb,m)$.
This means that $p|m$ and $p|a+tb$. Now suppose $p|a$, then $p|tb$ and by definition $p\not| t$ so $p|b$, but then
$p|Gcd(a,b,m)=1$, Contradiction. Now suppose $p\not| a$, then by definition $p|t$, so $p$ does divide $a$. Contradiction.\\
So there does not exist a prime $p$, that divides $Gcd(a+tb,m)$, which hence must be one.\\

Now the general case, $Gcd(a,b,m)=d$.
Define $a'=a/d$, $b'=b/d$ and $m'=m/d$, then $Gcd(a',b',m')=1$.
By the previous argument, there exists an $t$ such that $Gcd(a'+tb',m')=1$. Thus $Gcd(a+tb,m)=d$.
\end{proof}

\begin{lemma}\label{lemgen}
 Let $a,m\in \Z$, with $Gcd(a,m)=d$, then $\bar{a}$ is a generator of the additive subgroup $d\Z /m\Z$ of $\Z/m\Z$.
\end{lemma}
\begin{proof}
According to the Extended Euclidean Algorithm, there exist $u,v\in \Z$ such that $ua+vm=d$, so $\bar{u}\bar{a}=\bar{d}$.
Thus $\bar{u}\bar{a}$ is a generator of $d\Z/m\Z$, hence so is $\bar{a}$.\end{proof}

Now we will first prove a special case of Theorem \ref{main}, namely the three dimensional one, before we give the proof of the general case.
(The proof of this proposition is perhaps the most technical part of this article.)

\begin{proposition}\label{prop3var}
 Let $q=p^r$ and $\F_q$ the finite field with $q$ elements.
  \[\pi_{q^m}(\TA_3(\F _q))=\pi_{q^m}(\DA_3(\F_q))\]
\end{proposition}
\begin{proof}
It suffices to prove that $\pi_{q^m}((X+Y^aZ^b,Y,Z))\in\pi_{q^m}(\DA_3(\F_q))$ for $(a,b)\in \left(\mathbb{Z}/(q^m-1)\mathbb{Z}\right)^2$.
Let $\X = \F _{q^m}^3$ and $\tilde{\X}=X\setminus \{u\in \F_q^3|u_3=0\}$, then define the bijection $\psi:\tilde{\X}\rightarrow \tilde{\X}$ by $\psi (u_1,u_2,u_3)=(u_1u_3^{-1},u_2u_3,u_3)$.
Now from Proposition \ref{prop0} it follows that there exists an automorphisms $T_m\in DA_3(\F_{q})$, such that $\pi_{q^m}(T_m)|_{\tilde{\X}}=\psi$.

Suppose $\pi_{q^m}(E_{1,(\alpha,\beta)})\in\pi_{q^m}(\DA_3(\F_q))$, then
\begin{eqnarray}
\pi_{q^m}(T_m^{-1}E_{1,(\alpha,\beta)}T_m)&=&\pi_{q^m}(E_{1,(\alpha,\beta+\alpha +1)})\label{eqn1}\\
\pi_{q^m}(R_{2,3}E_{1,(\alpha,\beta)}R_{2,3})&=&\pi_{q^m}(E_{1,(\beta,\alpha)})\label{eqn2}
\end{eqnarray}
Equation (\ref{eqn1}) follows from Proposition \ref{prop1}.
Since $E_{1,(p-1,p-1)}\in \pi_{q^m}(\DA_3(\F_q))$ by definition, we need to prove that, starting with $E_{1,(p-1,p-1)}$ and applying (\ref{eqn1}) and (\ref{eqn2}), we can get $\pi_{q^m}(E_{1,(\alpha,\beta)})\in\pi_{q^m}(\DA_3(\F_q))$, for any pair $(\alpha, \beta)$.
The equations (\ref{eqn1}) and (\ref{eqn2}) translate into operations
\begin{eqnarray*}
 \varrho(\alpha,\beta) &=& (\alpha,\beta+\alpha+1)\\
\tau(\alpha,\beta) &=& (\beta,\alpha)
\end{eqnarray*}
on the space $\left(\mathbb{Z}/(q^m-1)\mathbb{Z}\right)^2$, where we can compute mod $(q^m-1)\mathbb{Z}$ as $\alpha^{q^m}=\alpha$ for any $\alpha \in \F_{q^m}$. So rephrasing the quoestion: Starting with $(\overline{p-1},\overline{p-1})\in \left(\mathbb{Z}/(q^m-1)\mathbb{Z}\right)^2)$ and iterating these two operation $\varrho$ and $\tau$, do we reach all of $\left(\mathbb{Z}/(q^m-1)\mathbb{Z}\right)^2$?

Unfortunately the answer is no. But we can reach almost every point and thereafter we show that we can mimic the maps we can not reach this way as well:
Suppose $Gcd(a+1,b+1,q^m-1)=1$, then from Lemma \ref{lemd} it follows that there exists a $t$ such that $Gcd(a+1+t(b+1),q^m-1)=1$.
From Lemma \ref{lemgen} and the fact that $Gcd(p,q^m-1)=1$ it follows that $p$ is a generator of the additive group $\mathbb{Z}/(q^m-1)\mathbb{Z}$ so there exists a $k_1$ such that $p-1+k_1p=a+t(b+1)$ mod $q^m-1$.
So from $\varrho^{k_1}((p-1,p-1))=(p-1,a+t(b+1))$ and since $Gcd(a+1+t(b+1),q^m-1)=1$, it follows that there exists a $k_2$ such that $p-1+k_2(a+1+t(b+1))=b$, thus $\varrho^{k_2}\tau\varrho^{k_1}((p-1,p-1))=(a+t(b+1),b)$
So $\tau\varrho^{q^m-1-t}\tau\varrho^{k_2}\tau\varrho^{k_1}((p-1,p-1))=(a,b)$.

Now if $Gcd(a+1,b+1,q^m-1)=d$, we need a little more work. From Lemma \ref{lemd} it follows that there exists a $t$ such that 
$Gcd(a+1+t(b+1),q^m-1)=d$. Now suppose we start with $(d-1,d-1)$, then we can write $a+t(b+1)=k_1d+(d-1)$ so $\varrho^{k_1}((d-1,d-1))=(d-1, a+t(b+1))$. Since $Gcd(a+1+t(b+1),q^m-1)=d$ it follows from Lemma \ref{lemgen} that $a+1+t(b+1)$ is a generator for the additive group $d \mathbb{Z}/(q^m-1)\mathbb{Z}$. Since $d|b+1$ it follows that there exists a $k_2$ such that $b=d-1+k_2(a+1+t(b+1))$, so $\varrho^{k_2}\tau\varrho^{k_1}((d-1,d-1))=(a+t(b+1),b)$.
So $\tau\varrho^{q^m-1-t}\tau\varrho^{k_2}\tau\varrho^{k_1}((d-1,d-1))=(a,b)$.

It remains to prove that we can reach $(d-1,d-1)$. Unfortunately this can not be done just using $\varrho$, $\tau$.
So we have to show that $\pi_{q^m}(E_{1,(d-1,d-1)})\in \pi_{q^m}(\DA_3(\F_q))$, where $d|q^m-1$. Now from the previous part it follows
that $\pi_{q^m}(E_{1,(d-1,d)})=\pi_{q^m}((X+Y^{d-1}Z^d,Y,Z))\in \pi_{q^m}(\DA_3(\F_q))$, since $Gcd(d,d+1,q^m-1)=1$.
Now let $d_1$ be the smallest divisor of $q^m-1$, then $\pi_{q^m}((X+(p-1)Y^{d_1-1}Z^{d_1},Y,Z))\pi_{q^m}((X,Y,Z+p-1))\pi_{q^m}((X+Y^{d_1-1}Z^{d_1},Y,Z))\pi_{q^m}((X,Y,Z+1))=
 \pi_{q^m}((X+(p-1)Y^{d_1-1}Z^{d_1} + Y^{d_1-1}(Z+1)^{d_1},Y,Z))=\pi_{q^m}((X+Y^{d_1-1}({d_1 \choose 1}Z^{d_1-1} +P(Z)),Y,Z)) $, with $\deg(P(Z))\leq d_1-2$. Since $d_1$ is the smallest divisor of $q^m-1$, it follows that $p$ does not divide $d_1$ and that $Gcd(d_1,d_1-i,q^m-1)=1$ for $1\leq i \leq d_1$.
From which it follows that $\pi_{q^m}((X+Y^{d_1-1}P(Z),Y,Z))\in  \pi_{q^m}(\DA_3(\F_q))$. Thus
\[ 
\begin{array}{rl}
&\pi_{q^m}\Big(
S_{1,d_1}
\big(X+Y^{d_1-1}P(Z),Y,Z\big) \big(X+Y^{d_1-1}({d_1 \choose 1}Z^{d_1-1} +P(Z)),Y,Z\big)\cdot S_{1,d_1^{-1}} \Big)\\
=&\pi_{q^m}((X+Y^{d_1-1}Z^{d_1-1},Y,Z)).
\end{array} \] Now let $d_2$ be the second smallest divisor, we can repeat the procedure described above and since we have already made all smaller degrees we have that $\pi_{q^m}((X+Y^{d_2-1}Z^{d_2-1},Y,Z))\in \pi_{q^m}(\DA_3(\F_q))$. Now by induction we are done.\end{proof}

Now we are ready to prove our main result, Theorem \ref{main}.\\
\begin{proof}
First of all note that it suffices to prove that $\pi_{q^m}(E_{1,\alpha})\in \pi_{q^m}(\DA_n(\F_q))$ for all $ \alpha \in \{0,\ldots ,q^m-1\}^{n-1}$.

Let $\X = \F _{q^m}^n$ and $\X_i=\X\setminus \{u\in \F_q^n|u_i=0\}$, then define the bijection $\psi_i:\X_i\rightarrow  \X_i$ by $\psi _i(u_1,\ldots ,u_n)=(u_1u_i^{-1},u_2u_i,u_3,\ldots ,u_n)$. Now Proposition \ref{prop0} states that there exist a map $T_m\in \DA_n(\F_q)$, such that
$\pi_{q^m}(T_m)|_{\X_n}=\psi_n $.
Now define $T_{m,i}=R_{i,n}T_mR_{i,n}$. Since $\pi_{q^m}$ is a group homomorphism, it follows that $\pi_{q^m}(T_{m,i})|_{\X_i}=\pi_{q^m}(R_{i,n}T_mR_{i,n})|_{\X_i}=\pi_{q^m}(R_{i,n})|_{\X_i}\pi_{q^m}(T_m)|_{\X_i}\pi_{q^m}(R_{i,n})|_{\X_i}=\psi_i$.\\
Now let $\alpha =(\alpha _2,\ldots ,\alpha _n) \in \{0,\ldots ,q^m-1\}^{n-1}$, for $i=4,\ldots ,n $ write $\alpha_i =pk_i+r_i$, with $r_i \in\{0\ldots ,p-1\}$. Then Lemma \ref{lem1} states that $\pi_{q^m}(E_{1,(p-1,p-1,r_4,\ldots ,r_n)})\in \pi_{q^m}(\DA_n(\F_q))$.
Now from Proposition \ref{prop1} it follows that $\pi_{q^m}(T_{m,i}^{-1}E_{1,(p-1,p-1,r_4,\ldots ,r_n)}T_{m,i})=\pi_{q^m}(E_{1,(p-1,p-1,r_4,\ldots ,r_i+p,\ldots ,r_n)})$.
So $\pi_{q^m}(T_{m,4}^{-k_4}\cdots T_{m,n}^{-k_n}E_{1,(p-1,p-1,r_4,\ldots ,r_n)}T_{m,n}^{k_n}\cdots T_{m,4}^{k_4})=\pi_{q^m}(E_{1,(p-1,p-1,\alpha_4,\ldots ,\alpha_n)})\in \pi_{q^m}(\DA_n(\F_q)$.
To prove the final step one can copy the proof of Proposition \ref{prop3var}, and extend all automorphisms with $n-3$ variables.\end{proof}

\section{Tamizables versus Linearizables}
\label{lin}

In this section we will compare two subgroups of $\GA_n(k)$ ($k$ a field), namely the group generated by the so called linearizables ($\GLIN_n(k)$) and the group generated by the tamizables ($\GTAM_n(k)$), both introduced in \cite{PolMau}. An automorphism $F\in \GA_n(k)$ is called linearizable if it is the conjugate of a linear automorphism, so if there exist an $L\in \GL_n(k)$ and a $G\in \GA(k)$, such that $F=G^{-1}\circ L\circ G$.
Similarly, an automorphism is called tamizable if it is the conjugate of a tame automorphism.
\begin{definition}
 Let $G$ be a group, and $H$ a subgroup of $G$. We define $\NN (H,G)$ to be the smallest normal subgroup of $G$ that contains $H$, i.e.
\begin{displaymath}
 \NN (H,G) = \langle g^{-1}hg|h\in H, g \in G \rangle
\end{displaymath}
\end{definition}
Furthermore let $g,h \in G$ then we write the commutator as $[ g,h] :=g^{-1}h^{-1}gh $.

Now we can define the following subgroups of $\GA_n(k)$:
\begin{eqnarray*}
\GLIN_n(k) & := & \NN (\GL_n(k),\GA_n(k))\\
\GTAM_n(k) & := & \NN (\TA_n(k),\GA_n(k))\\
\TLIN_n(k) & :=& \NN(\GL(k), \TA_n(k)\\ 
\end{eqnarray*}
(Note that some ``$\textup{TTAM}$'' would equal $\TA_n(k)$.)
Then the following is obviously true:\\
\ \\
\begin{tabular}{cccccc}
& $\TA_n(k)$ & $\subseteq $ & $ \GTAM_n(k) $ &&\\
 && & &\begin{turn}{-25}\parbox[t]{5mm}{$\subseteq $ }\end{turn} & \\
&\begin{turn}{90}\parbox[t]{5mm}{$\subseteq $ }\end{turn} & &\begin{turn}{90}\parbox[t]{5mm}{$\subseteq $ }\end{turn} & & $\GA_n(k)$\\ 
 && & &\begin{turn}{25}\parbox[b]{5mm}{$\subseteq $ }\end{turn} & \\
 $\GL_n(k)$ $\subseteq $ &$\TLIN_n(k)$ & $\subseteq $ & $ \GLIN_n(k) $ &&\\
\end{tabular}\ \\ \ \\

One of the  motivations of this section is in attacking the so-called {\em linearization problem}, which is the conjecture that if $F^s=I$, then $F$ is linearizable. In particular, $F\in \GLIN_n(k)$. 
Note that in characteristic $p$ one actually has non-linearizable automorphisms like $F:=(X+Y^2,Y)$, for which $F^p=I$ and indeed $F$ is non-linearizable. (Or many an automorphism of an additive group action, for that matter.) This indicates that the linearization problem over characteristic $p$ should be reformulated:\\
\ \\
{\bf Linearization problem} in characteristic $p$: let $F\in \GA_n(k)$ where $k$ is a field of characteristic zero. Assume that $F^s=I$ where $Gcd(p,s)=1$. Then $F$ is linearizable. \\

As observed above, an automorphism $F^s=I$, $Gcd(p,s)=1$ for which $F\not \in \GLIN_n(k)$, must be a counterexample to the above problem.
However, in this section we show that such an approach cannot work except if $k=\F_2$: if $k\not = \F_2$, then $\GLIN_n(k)=\GTAM_n(k)$, the latter being a good candidate of equalling $\GA_n(k)$. 
The main result of this section is the following: 

\begin{theorem}\label{LinTamMain}
If $n\geq 1$ and $k\not = \F_2$, then $\GLIN_n(k)=\GTAM_n(k)$ and $\TLIN_n(k)=\TA_n(k)$.
In case $k=\F_2$, then $\GLIN_n(\F_2)\subsetneq \GTAM_n(\F_2)$ and $\TLIN_n(\F_2)\subsetneq \TA_n(\F_2)$ 
\end{theorem}

\begin{proof}
The theorem follows directly from lemma \ref{lemma2} and proposition \ref{PropB}, and the implication 
``$\TLIN_n(k)=\TA_n(k)$ then $\GLIN_n(k)=\GTAM_n(k)$'' and its negation ``$\GLIN_n(k)\subsetneq \GTAM_n(k)$ then $\TLIN_n(k)\subsetneq \TAM_n(k)$''.
\end{proof}

%\subsection{The case $k\not= \F_2$}\label{qnot2}

\begin{lemma}\label{lemma2}
If $k\not = \F_2$, then 
 $\TA_n(k)= \TLIN_n(k)$.
\end{lemma}
\begin{proof}
Only the inclusion ``$\subsetneq$'' needs to be proven.
Since $\GL_n(k)<\TLIN_n(k)$ it suffices to prove that $E_{1,\alpha}\in \TLIN_n(k)$ for all $\alpha \in \N^{n-1}$.
Choose $c\not = 0,1$ (which is possible since $k\not = \F_2$) and $d:=(1-c)^{-1}$.
Then an elementary computation shows that
\begin{eqnarray*}
S_{1,c}\left(S_{1,d}E_{1,\alpha}S_{1,d^{-1}}\right)^{-1}S_{1,c}^{-1}\left(S_{1,d}E_{1,\alpha}S_{1,d^{-1}}\right) & = & E_{1,\alpha}
\end{eqnarray*}
Since $S_{1,c}$ and $\left(S_{1,d}E_{1,\alpha}S_{1,d^{-1}}\right)^{-1}S_{1,c}^{-1}\left(S_{1,d}E_{1,\alpha}S_{1,d^{-1}}\right)$ are in $\TLIN_n(k)$ by definition, it follows that $E_{1,\alpha} \in \TLIN_n(k)$.
\end{proof}

%\subsection{The case $k=\F_2$ }\label{qis2}
%In this section we will assume that $k=\F_2$. The aim of this section is to prove:
\begin{proposition}\label{PropB}
 $\GTAM_n(\FF _2) \not= \GLIN_n(\FF _2),\ n\geq 2$.
\end{proposition}

\begin{proof} This follows from Lemma \ref{subset} and the appendix.
\end{proof}

%\subsubsection{$n\geq3$}\label{qis2ngeq3}
%In case $n\geq3$ we will show that $\GTAM_n(\FF _2) \not= \GLIN_n(\FF _2)$.\\
Let $\mathcal{A}_{2^n}$ be the alternating subgroup of the symmetric group $\mathcal{S}_{2^n}\cong \Bij(\F_2^n)$. 
\begin{lemma}\label{subset}
 $\pi_2(\GLIN_n(\FF _2)) \subseteq \mathcal{A}_{2^n}$.
\end{lemma}
\begin{proof}
First remark that if $\pi_2(h)\in \mathcal{A} _{2^n}$, then also $\pi_2(g^{-1}hg)\in \mathcal{A} _{2^n}$, because $\pi_2(gh)=\pi_2(g)\pi_2(h)$.
So all we have to do is show that $\pi_2(\GL_n(\FF _2))\subseteq \mathcal{A}_{2^n}$. 
Because $\mathcal{A}_{2^n}$ is a group it suffices to prove that for the generators of $\GL_n(\FF _2)$:
\begin{eqnarray*}
F_1 & := & (x_1+x_2,x_3,\ldots ,x_n)\\
F_i & := & R_{1,i} = (x_i,x_2,\ldots,x_{i-1},x_1,x_{i+1},\ldots ,x_n)\ for\ 2\leq i\leq n 
\end{eqnarray*}
$\pi_2(F_i)\in \mathcal{A}_{2^n}$, $\mathcal{a} \in \mathcal{A}_{2^n} \Leftrightarrow \mathcal{a}$ is the compostion of an even number of 2 cycles.
The number of 2 cycles that appear in $\pi_2(F_i)$ is equal to $\frac{2^n-\#Fix(F_i)}{2}$ where $2^n$ is the number of elements in $\FF _2^n$ and $Fix(F_i)$ is the set of fixed points in $\FF _2^n$ under $F_i$. But $\#Fix(F_i)=2^{n-1}$, so the number of 2 cycles in $F_i$ is $\frac{2^n-2^{n-1}}{2} = 2^{n-2}$ which is even for $n\geq 3$. So $F_i\in \mathcal{A}_{2^n}\forall i$. Which proves the claim.
\end{proof} 

\begin{proposition}
 $\GTAM_n(\FF _2) \not= \GLIN_n(\FF _2)$ if $n\geq 3$.
\end{proposition}

\begin{proof}
Theorem 2.3 in \cite{Maubach1} states that $\pi_2 (\TA_n(\F_2)) = \Bij(\F_2^n)$.
So $\Bij(\F_2^n) = \pi_2 (\TA_n(F_2))\subseteq \pi_2(\GTAM_n(\FF_2)) \subseteq \Bij(\F_2^n)$ it follows that $\pi_2(\GTAM_n(\FF_2)) = \Bij(\F_2^n)$.
Furthermore we have just shown in Lemma \ref{subset} that  $\pi_2(\GLIN_n(\FF _2)) \subseteq \mathcal{A}_{2^n}$, but $\mathcal{A}_{2^n}\not= \Bij(\F_2^n)$.
From which we can conclude that  $\GTAM_n(\FF _2) \not= \GLIN_n(\FF _2)$.
\end{proof}

%\subsubsection{$n=2$}\label{qis2nis2}
%First recall that not only does an element of $\GA_n(\FF _q)$ induce a bijection of $\FF_q^n$, 
%it also induces a bijection of $\FF_{q^m}^n$ for $m\geq 2$. Now we are going to look at the special case where $q=2$ and $n=2$.
%So we can look at the functor:
%\begin{displaymath}
%\pi _4:\GA_2(\F_2)\rightarrow \mathcal{B}(\F_4^2),
%\end{displaymath}
%as the functor sending $F\in GA_2(\FF_2)$ to the map $\pi_4(F):\F_4^2\rightarrow \FF_4^2$.
%We will show that
\begin{proposition}\label{q2n2}
 $\GTAM_2(\FF _2) \not= \GLIN_2(\FF _2)$
\end{proposition}
\begin{proof} 
It follows from Jung-van der Kulk, that $\TA_2(\FF_2)=\GTAM_2(\FF _2)=\GA_2(\FF_2)$. Furthermore $\pi_4$ is a grouphomomorphism, so $\pi_4$ commutes with $\NN$.
We have shown using MAGMA that $[\pi_4(\GLIN_2(\FF_2)):\pi_4(\GA_2(\FF_2))]=2$,  the program we used to compute this can be found below.
\end{proof}

Here follows the MAGMA program we used to prove that\\ $[\pi_4(\GLIN_2(\FF_2)):\pi_4(\GA_2(\FF_2))]=2$.\\
\begin{verbatim}
/*Create a finite field with 4 elements.*/
F4:=GF(4);

/*Create a polynomialRing in two variables over F4.*/
P<x,y>:=PolynomialRing(F4,2);

/*Create an ordered set with all elements of F4^2.*/
J:=\{@ [a,b]:a,b\ in F4 @};

/*Create a list of all generating polynomials of the tame  
automorphism group over F4, with coefficients in F2. Note 
that since we are only interested in the bijections they 
induce on F4^2, it suffices to give generators up to degree 3.*/
GeneratingPolynomialsTA:=[[x+1,y],[x+y,y],[x+y^2,y],[x+y^3,y],[y,x]];

/*Create a list of the generating polynomials of the linear 
group over F4, with coefficients in F2.*/
GeneratingPolynomialsGL:=[[x+y,y],[y,x]];

/*Create the subgroup of the symmetricgroup of F4^2, 
generated by bijections induced by the tame automorphisms, 
with coefficients in F2.*/
Pi2TA:=sub<Sym(4^2)|[[Index(J,[Evaluate(P[1],[a,b]),Evaluate(P[2],
[a,b])]):a,b in F4]:P in GeneratingPolynomialsTA]>;

/*Create the subgroup of Pi2TA, generated by the bijections 
induced by linear automorphisms with coefficients in F2.*/
Pi2GL:=sub<Pi2TA|[[Index(J,[Evaluate(P[1],[a,b]),Evaluate(P[2],
[a,b])]):a,b in F4]:P in GeneratingPolynomialsGL]>;

/*Create the normal closure (or conjugate closure) of Pi2GL in 
Pi2TA.*/
NCPi2GL:=NormalClosure(Pi2TA,Pi2GL);

/*Show that the #Pi2TA / #NCPi2GL = 2*/
#Pi2TA / #NCPi2GL; 
\end{verbatim}

\section{Open problems}

We gather two  of the open problems we found on the way:

\begin{problem} Does there exist a finite set $E\subset \TA_n(\F_q)$ such that $<\Aff_n(\F_q),E>=\TA_n(\F_q)$? ($n\geq 3$)
\end{problem}

\begin{problem}
Is $\DA_n(\F_q)$ indeed non-equal to $\TA_n(\F_q)$? ($n\geq 3$)
\end{problem}

\end{document}